\documentclass[]{article}
\begin{document}
\newtheorem{proposition}{Proposition}[section]
\newtheorem{definition}{Definition}[section]
\newtheorem{lemma}{Lemma}[section]

\title{\bf Generalization of One of Lie's Theorems}
\author{Keqin Liu\\Department of Mathematics\\The University of British Columbia\\Vancouver, BC\\
Canada, V6T 1Z2}
\date{January 31, 2008}
\maketitle

\begin{abstract} We prove a generalization of one of Lie's Theorems in the context of 
Lie-like algebras$^{2-nd}$.\end{abstract}

The theory of Lie algebras has many applications in mathematics and physics. One possible way of generalizing the theory of Lie algebras is to develop the theory of Lie-like algebras$^{2-nd}$ algebras, where the notion of a Lie-like algebras$^{2-nd}$ algebra was introduced in \cite{Liu14}.
One of Lie's Theorems claims that the only irreducible representations of a solvable Lie algebra over an algebraically closed field $\mathbf{k}$ of characteristic $0$ have dimension $1$. We call this theorem 
Lie's Theorem for convenience in this paper.
Since Lie's Theorem  is one of the fundamental results in the theory of Lie algebras, finding the counterpart of Lie's Theorem in the context of Lie-like algebras$^{2-nd}$ is of importance to develop the theory of Lie-like algebras$^{2-nd}$. The purpose of this paper is to prove a generalization of Lie's Theorem in the context of Lie-like algebras$^{2-nd}$.

\medskip
In this paper, all vector spaces are vector spaces over an algebraically closed field $\mathbf{k}$ of characteristic $0$. The generalization of Lie's Theorem in this paper claims that a finite dimensional nonzero ordinary module over a finite dimensional solvable Lie-like algebras$^{2-nd}$ always contains a one-dimensional ordinary submodule. Since a finite dimensional nonzero module over a finite dimensional solvable Lie superalgebra does not always contain a one-dimensional submodule (\cite{Kac}), Lie-like algebras$^{2-nd}$ are closer to Lie algebras than Lie superalgebras.  Hence, using Lie-like algebras$^{2-nd}$ seems to be more suitable than using Lie superalgebras in order to extend Lie theory. 

\medskip
After discussing the basic properties of modules over Lie-like algebras$^{2-nd}$ in Section 1, we use Section 2 to give the proof of the generalization of Lie's Theorem.

\medskip
\section{Modules over Lie-like Algebras$^{2-nd}$}

We begin this section with the definitions of a Lie-like algebra$^{2-nd}$ $L$ and an ordinary 
$L$-module(\cite{Liu14}).

\begin{definition}\label{def1.1} Let $S$ be a nonempty set. A vector space $L$ is called a 
{\bf Lie-like algebra$^{2-nd}$} induced by the set $S$ if there exist a family of binary operations
$$
 \Big\{\, \langle\; , \;\rangle _k \;\Big |\;\mbox{$\langle\; , \;\rangle _k : L\times L\to L$ is a binary map and $k\in S$} \,\Big\}
$$
such that both the {\bf Jacobi-like identity$^{2-nd}$}
\begin{equation}\label{eq1.1}
\langle \langle x, y\rangle _k, z\rangle _h=\langle x, \langle y, z\rangle _h\rangle _k+
\langle \langle x, z\rangle _h, y\rangle _k
\end{equation}
and the following identity
\begin{equation}\label{eq1.1'}
\langle \langle x, y\rangle _k, z\rangle _h=\langle \langle x, y\rangle _h, z\rangle _k
\end{equation}
hold for $x$, $y$, $z\in L$ and $h$, $k\in S$
\end{definition}

A Lie-like algebra$^{2-nd}$ $L$  induced by the set $S$ is also denoted by 
$(\, L\, , \langle \, , \, \rangle_{k\in S}\,)$. The notion of a Leibniz algebra was introduced in \cite{Loday}. Clearly, if $L$ is a Lie-like algebra$^{2-nd}$ induced by a set $S$, then $L$ is a Leibniz algebra with respect to the binary operation $\langle\; , \;\rangle _k$ for each $k\in S$. Thus, a Leibniz algebra is a Lie-like algebra$^{2-nd}$ $L$  induced by the set whose cardinal number is $1$, and a Lie-like algebra$^{2-nd}$ induced by a set $S$ is a bundle of Leibniz algebras satisfying the Jacobi-like identity$^{2-nd}$ and (\ref{eq1.1'}). 

\medskip
A Lie-like algebra$^{2-nd}$ $(\, L\, , \langle \, , \, \rangle_{k\in S}\,)$ is said to be
{\bf non-trivial} if there do not exist a binary operation $\langle \, , \, \rangle: L\times L\to L$ and a map $\phi : S\to {\bf k}$ such that
\begin{equation}\label{eq1.2}
 \langle x, y\rangle _k=\phi (k)\langle x, y\rangle\quad\mbox{for $x$, $y\in L$ and $k\in S$}.
\end{equation}
A subspace $I$ of a Lie-like algebra$^{2-nd}$ $(\, L\, , \langle \, , \, \rangle_{k\in S}\,)$
is called an {\bf ideal} of $L$ if $\langle I, L\rangle_s\subseteq I$ and 
$\langle L, I\rangle_s\subseteq I$ for all $s\in S$.

\begin{definition}\label{def1.5}  Let $L$ be a Lie-like algebra$^{2-nd}$ induced by a set $S$. A vector space $V$ is called an {\bf ordinary module} over $L$ (or an {\bf ordinary $L$-module}) if there exist a family of linear maps
$$
 \bigg\{\, f_k, \; g_k \;\bigg |\; \begin{array}{l}\mbox{$f_k : x\mapsto f_k(x)$ and 
$g_k : x\mapsto g_k(x)$ are linear}\\\mbox{maps from $L$ to $End(V)$ for $x\in L$ and $k\in S$}
\end{array} \,\bigg\}
$$
such that
\begin{equation}\label{eq1.3}
 f_h(\langle x, y\rangle _k)=[f_h(x), \, f_k(y)],
\end{equation}
\begin{equation}\label{eq1.4}
 g_h(\langle x, y\rangle _k)=[g_h(x),\, f_k(y)],
\end{equation}
\begin{equation}\label{eq1.5}
g_k(x)g_h(y)=g_h(x)f_k(y)=g_k(x)f_h(y),
\end{equation}
and
\begin{equation}\label{eq1.6}
f_k(x)f_h(y)=f_h(x)f_k(y),\quad f_k(x)g_h(y)=f_h(x)g_k(y),
\end{equation}
where $x$, $y\in L$, and $[f_h(x), f_k(y)]:=f_h(x)f_k(y)-f_k(y)f_h(x)$ is the ordinary bracket. 
An ordinary $L$-module $V$ is also denoted by $(V, \{f_k\}, \{g_k\}_{k\in S})$, and $f_k$ and 
$g_k$ are called the {\bf right linear map indexed by $k$} and the {\bf left linear map indexed by $k$} of $V$, respectively.
\end{definition}

The notion of a module over a Leibniz algebra, which was introduced in \cite{LT}, is a special case of the notion of an ordinary module over a Lie-like algebra$^{2-nd}$. Clearly, an ordinary module over a Lie-like algebra$^{2-nd}$ is a bundle of modules over Leibniz algebras.

\medskip
If $(V, \{f_k\}, \{g_k\}_{k\in S})$ is an ordinary module over a Lie-like algebra$^{2-nd}$ $L$, then
\begin{equation}\label{eq1.7a}
f_h(\langle x, y\rangle _k)=f_k(\langle x, y\rangle _h)
\quad\mbox{for $x$, $y\in L$ and $h$, $k\in S$}
\end{equation}
and
\begin{equation}\label{eq1.7b}
g_h(\langle x, y\rangle _k)=g_k(\langle x, y\rangle _h)
\quad\mbox{for $x$, $y\in L$ and $h$, $k\in S$.}
\end{equation}

\medskip
If $(\, L\, , \langle \, , \, \rangle_{k\in S}\,)$ is a Lie-like algebra$^{2-nd}$, then 
$(L, -r_k, \ell_k)$ is an ordinary module over the  Lie-like algebra$^{2-nd}$ $L$, where $r_k$ and $\ell_k$ are the right multiplication and the left multiplication of $L$, respectively; that is,
\begin{equation}\label{eq1.7}
r_k(x)(a): =\langle a , \, x \rangle_k\quad\mbox{and}\quad 
\ell_k(x)(a): =\langle x , \, a \rangle_k\quad\mbox{for $x$, $a\in L$ and $k\in S$.}
\end{equation}
This ordinary $L$-module $(L, -r_k, \ell_k)$ is called the {\bf adjoint module} over the Lie-like algebra$^{2-nd}$ $L$.

\medskip
Let $(V, \{f_k\}, \{g_k\}_{k\in S})$ be an ordinary module over a Lie-like algebra$^{2-nd}$
$(\, L\, , \langle \, , \, \rangle_{k\in S}\,)$. A subspace $U$ of $V$ is called an 
{\bf ordinary submodule} of $V$ if 
$$
f_k(x)(U)\subseteq U\quad\mbox{and}\quad g_k(x)(U)\subseteq U\quad\mbox{for 
$x\in L$ and $k\in S$}. 
$$
Following \cite{Liu14}, the subspace
\begin{equation}\label{eq1.9}
V^{ann,+}_{2-nd}:=\sum_{x\in L,\; v\in V\atop h, k\in S}{\bf k}(g_h(x)-f_k(x))(v)
\end{equation}
is called the {\bf plus annihilator} of $V$, which is an ordinary submodule of $V$ and plays an important role in the proof of the generalization of Lie's Theorem. 

\medskip
A Lie-like algebra$^{2-nd}$
$(\, L\, , \langle \, , \, \rangle_{k\in S}\,)$ is said to be {\bf solvable} if there exists a positive integer $n$ such that $\mathcal{D}^n\mathcal{L}=0$, where $\mathcal{D}^n\mathcal{L}$ is defined inductively as follows:
\begin{equation}\label{eq1.10}
\mathcal{D}^1\mathcal{L}:=\mathcal{L}, \quad 
\mathcal{D}^{n+1}\mathcal{L}=\sum_{k\in S}\langle\mathcal{D}^n\mathcal{L}, \,\mathcal{D}^n\mathcal{L}\rangle_k\quad\forall n\ge 1.
\end{equation}

\medskip
\section{The Proof of the Main Theorem}

In order to present the proof of the generalization of Lie's Theorem in a clear way, we divide it into a few parts. First, we have

\begin{proposition}\label{pr2.1} Let $V$ be a finite dimensional vector space over $\mathbf{k}$, and let $\mathcal{G}$ be a subspace of the general linear Lie algebra $g\ell(V)$. If $\mathcal{A}$ is a subspace of $g\ell(V)$ and
$$ \mathcal{G}\subseteq \mathcal{N}_{g\ell(V)}(\mathcal{A}):=\,\{\, X\in g\ell(V)\, |\, [\mathcal{A}, X]\subseteq \mathcal{A} \, \},
$$
then
$$
U:=\,\{\, u\in V \, |\, A(u)=\phi (A)u \quad\mbox{for all $A\in\mathcal{A}$} \, \}
$$
is an invariant subspace of $\mathcal{G}$; that is, $\mathcal{G}(U)\subseteq U$, where 
$\phi : \mathcal{A}\to \mathbf{k}$ is a linear functional.
\end{proposition}

\medskip
\noindent
{\bf Proof} The proof of Proposition 2.1 is the same as the proof of Lemma 2 on page 24 in \cite{Wan}.
 
\hfill\raisebox{1mm}{\framebox[2mm]{}}

\bigskip
Next, we have

\begin{proposition}\label{pr2.2} Let $L=\mathcal{A}\oplus \mathbf{k}x$ (direct sum of vector spaces) be a Lie-like algebra$^{2-nd}$ over the field $\mathbf{k}$, where $\mathcal{A}$ is an ideal of $L$ and 
$\langle L, L\rangle_k\subseteq \mathcal{A}$ for $k\in S$. Let $(V, \{f_k\}, \{g_k\}_{k\in S})$  be a finite dimensional ordinary $L$-module such that
\begin{equation}\label{eq2.6}
f_k(A)(u_0)=\phi_k (A) u_0, \quad g_k(A)(u_0)=\psi_k (A) u_0\quad\mbox{for $A\in\mathcal{A}$ 
and $k\in S$,}
\end{equation}
where $u_0\in V$, $\phi_k : \mathcal{A}\to \mathbf{k}$ and $\psi_k : \mathcal{A}\to \mathbf{k}$ are  functionals. If $h\in S$ and $u_m:=g_h(x)^m(u_0)$ for $m\in\mathcal{Z}_{\ge 0}$, then
\begin{eqnarray}
\label{eq2.7}f_k(A)(u_m)& \equiv & \phi_k (A)u_m \; mod \; 
\Big(V^{ann,+}_{2-nd}+\displaystyle\sum_{i=0}^{m-1}\mathbf{k}u_i\Big), \\
\label{eq2.8}g_k(A)(u_m)& \equiv & \psi_k (A)u_m \; mod \; 
\Big(V^{ann,+}_{2-nd}+\displaystyle\sum_{i=0}^{m-1}\mathbf{k}u_i\Big) 
\end{eqnarray}
and
\begin{equation}\label{eq2.9}
f_k(x)(u_m)\,\equiv\,u_{m+1} 
\; mod \; \Big(V^{ann,+}_{2-nd}+\displaystyle\sum_{i=0}^m\mathbf{k}u_i\Big)
\end{equation}
for all $A\in \mathcal{A}$ and  $k\in S$.
\end{proposition}

\medskip
\noindent
{\bf Proof} We use induction on $m$ to prove (\ref{eq2.7}), (\ref{eq2.8}) and (\ref{eq2.9}). If $m=0$, then  both (\ref{eq2.7}) and (\ref{eq2.8}) hold by  (\ref{eq2.6}), and (\ref{eq2.9}) holds by  (\ref{eq1.9}). 

\medskip
We now prove that if (\ref{eq2.7}), (\ref{eq2.8}) and (\ref{eq2.9}) hold for a fixed nonnegative integer $n$ with $n\le m$, then
\begin{eqnarray}
\label{eq2.10}f_k(A)(u_{m+1})& \equiv & \phi_k (A)u_{m+1} \, mod \, 
\Big(V^{ann,+}_{2-nd}+\displaystyle\sum_{i=0}^m\mathbf{k}u_i\Big),\\
\label{eq2.11}g_k(A)(u_{m+1})& \equiv & \psi_k (A)u_{m+1}\, mod \, 
\Big(V^{ann,+}_{2-nd}+\displaystyle\sum_{i=0}^m\mathbf{k}u_i\Big) 
\end{eqnarray}
and
\begin{equation}\label{eq2.12}
f_k(x)(u_{m+1})\,\equiv\,u_{m+2} 
\, mod \, \Big(V^{ann,+}_{2-nd}+\displaystyle\sum_{i=0}^{m+1}\mathbf{k}u_i\Big) 
\end{equation}
for all $A\in \mathcal{A}$ and  $k\in S$.

\medskip
Since $V^{ann,+}_{2-nd}$ is a submodule of $V$, we have
\begin{equation}\label{eq2.13}
f_k(z)(V^{ann,+}_{2-nd})\subseteq V^{ann,+}_{2-nd} \quad\mbox{for $z\in L$ and  $k\in S$.}
\end{equation}

\medskip
Using (\ref{eq1.4}), we have
\begin{eqnarray}
f_k(A)(u_{m+1})&=&f_k(A)g_h(x)(u_m) =\Big(-[g_h(x), f_k(A)]+g_h(x)f_k(A)\Big)(u_m)\nonumber\\
\label{eq2.14}&=&-g_h(\langle x, A\rangle_k )(u_m)+g_h(x)f_k(A)(u_m).
\end{eqnarray}

Since $\mathcal{A}$ is an ideal of $L$, $\langle x, A\rangle_k\in \mathcal{A}$. Hence, we get
$$
g_h(\langle x, A\rangle_k )(u_m)\stackrel{(\ref{eq2.8})}{\equiv}
\psi_h (\langle  x, A\rangle_k)u_m \; mod \; 
\Big(V^{ann,+}_{2-nd}+\displaystyle\sum_{i=0}^{m-1}\mathbf{k}u_i\Big)
$$
or 
\begin{equation}\label{eq2.15}
g_h(\langle x, A\rangle_k )(u_m)\,\equiv\,
0\; mod \; \Big(V^{ann,+}_{2-nd}+\displaystyle\sum_{i=0}^m\mathbf{k}u_i\Big).
\end{equation}

Using the fact that $g_h(x)(V^{ann,+}_{2-nd})\subseteq V^{ann,+}_{2-nd}$, we get
\begin{eqnarray}
&&g_h(x)\Big(f_k(A)(u_m)\Big)\nonumber\\
&\stackrel{(\ref{eq2.7})}{\equiv}
&g_h(x)\Big(\phi_k (A)u_m\Big)\;mod \;
\Big(g_h(x)(V^{ann,+}_{2-nd})+\displaystyle\sum_{i=0}^{m-1}\mathbf{k}g_h(x)(u_i)\Big) \nonumber\\
&\,\equiv\, &\phi_k (A)g_h(x)(u_m) \;mod \;
\Big(V^{ann,+}_{2-nd}+\displaystyle\sum_{i=1}^m\mathbf{k}u_i\Big) \nonumber\\
\label{eq2.16}&\,\equiv\, &\phi_k (A)u_{m+1} \;mod \;
\Big(V^{ann,+}_{2-nd}+\displaystyle\sum_{i=1}^m\mathbf{k}u_i\Big).
\end{eqnarray}

It follows from (\ref{eq2.14}), (\ref{eq2.15}) and (\ref{eq2.16}) that (\ref{eq2.10}) holds.

\medskip
Similarly, one can prove that (\ref{eq2.11}) holds.

\medskip
Finally, using the fact that $\langle x, x\rangle_k \in \mathcal{A}$, we have
\begin{eqnarray*}
f_k(x)(u_{m+1})&=&f_k(x)g_h(x)(u_m)=\Big(-[g_h(x), f_k(x)]+g_h(x)f_k(x)\Big)(u_m)\\
&\stackrel{(\ref{eq1.4})}{=}&-g_h(\langle x, x\rangle_k)(u_m)+g_h(x)f_k(x)(u_m)\\
&\stackrel{(\ref{eq2.8})}{\equiv} & 0+ g_h(x)f_k(x)(u_m)\; mod \; 
\Big(V^{ann,+}_{2-nd}+\displaystyle\sum_{i=0}^m\mathbf{k}u_i\Big) \\
&\stackrel{(\ref{eq2.9})}{\equiv} & g_h(x)(u_{m+1})\; mod \; 
\Big(g_h(x)(V^{ann,+}_{2-nd})+\displaystyle\sum_{i=0}^m\mathbf{k}g_h(x)(u_i)\Big) \\
&\equiv & u_{m+2}\; mod \; 
\Big(\displaystyle\sum_{i=1}^{m+1}\mathbf{k}u_i\Big),
\end{eqnarray*}
which proves that (\ref{eq2.12}) holds.

\hfill\raisebox{1mm}{\framebox[2mm]{}}

\bigskip
Using Proposition~\ref{pr2.2}, we now prove the following

\begin{proposition}\label{pr2.3} Let $L=\mathcal{A}\oplus \mathbf{k}x$ (direct sum of vector spaces) be a 
Lie-like algebra$^{2-nd}$ over the field $\mathbf{k}$, where $\mathcal{A}$ is an ideal of $L$ and 
$\langle L, L\rangle_k\subseteq \mathcal{A}$ for $k\in S$. Let $\phi_k: \mathcal{A}\to \mathbf{k}$ and 
$\psi_k : \mathcal{A}\to \mathbf{k}$ be functionals for $k\in S$. If $(V, \{f_k\}, \{g_k\}_{k\in S})$  is a finite dimensional ordinary $L$-module, then
\begin{equation}\label{eq2.21}
f_h(x)(U_{\phi}\cap U_{\psi})\subseteq U_{\phi}\cap U_{\psi}\quad\mbox{for all $h\in S$}
\end{equation}
and
\begin{equation}\label{eq2.22}
g_h(x)(u)\subseteq U_{\phi}
\quad\mbox{for all $h\in S$ and $u\in (U_{\phi}\cap U_{\psi})\setminus V^{ann,+}_{2-nd}$ },
\end{equation}
where 
$$
U_{\phi}:=\bigcap_{k\in S}U_{\phi_k},\quad\quad U_{\psi}:=\bigcap_{k\in S}U_{\psi_k},
$$
$$
U_{\phi_k}:=\{\,u\in V \, |\, \mbox{$f_k(A)(u)=\phi_k (A) u$ for $A\in \mathcal{A}$}  \,\}
$$
and
$$
U_{\psi_k}:=\{\,u\in V \, |\, \mbox{$g_k(A)(u)=\psi_k (A) u$ for $A\in \mathcal{A}$}  \,\}.
$$
\end{proposition}

\medskip
\noindent
{\bf Proof} Since $\mathcal{A}$ is an ideal of $L$, we have
\begin{equation}\label{eq2.23}
[f_h(x), f_k(A)]\stackrel{(\ref{eq1.3}) \& (\ref{eq1.7a})}{=}f_k(\langle x, A\rangle_h)\in f_k(\mathcal{A})\Rightarrow
f_h(x)\in \mathcal{N}_{g\ell(V)}(f_k(\mathcal{A}))
\end{equation}
and
\begin{equation}\label{eq2.24}
[g_k(A), f_h(x)]\stackrel{(\ref{eq1.4})}{=}g_k(\langle A, x\rangle_h)\in g_k(\mathcal{A})\Rightarrow
f_h(x)\in \mathcal{N}_{g\ell(V)}(g_k(\mathcal{A})).
\end{equation}

By (\ref{eq2.23}) and Proposition~\ref{pr2.1}, $f_h(x)$ does not change $U_{\phi_k}$, which implies that
\begin{equation}\label{eq2.25}
f_k(A)\Big(f_h(x)(u)\Big)=\phi_k (A) f_h(x)(u)\quad\mbox{for $A\in \mathcal{A}$, $u\in U_{\phi_k}$ and $k$, $h\in S$.}
\end{equation}

Similarly, by (\ref{eq2.24}) and Proposition~\ref{pr2.1}, $f_h(x)$ does not change $U_{\psi_k}$, which implies that
\begin{equation}\label{eq2.26}
g_k(A)\Big(f_h(x)(u)\Big)=\psi_k (A) f_h(x)(u)\quad\mbox{for $A\in \mathcal{A}$, $u\in U_{\psi_k}$  and $k$, $h\in S$.}
\end{equation}

Using (\ref{eq2.25}) and (\ref{eq2.26}), we get (\ref{eq2.21}).

\bigskip
We now prove
\begin{equation}\label{eq2.27}
f_k(A)\Big(g_h(x)(u)\Big)=\phi_k (A) g_h(x)(u)
\end{equation}
for $A\in \mathcal{A}$, $k$, $h\in S$ and 
$u\in (U_{\phi}\cap U_{\psi})\setminus V^{ann,+}_{2-nd}$.

Let $h$ be an arbitrary element of $S$. We define $u_0:=u$ and $u_m:=g_h(x)^m(u_0)$ for
$m\in\mathcal{Z}_{\ge 0}$. Since $V$ is finite dimensional and $u_m\in V$ for $m\in\mathcal{Z}_{\ge 0}$, there exists a positive integer $p$ such that
\begin{equation}\label{eq2.28}
\mbox{$W_h:=V^{ann,+}_{2-nd}\oplus {\bf k}u_0\oplus\cdots \oplus{\bf k} u_{p-1}$ and $u_m\in W$ for
$m\in\mathcal{Z}_{\ge 0}$.}
\end{equation}

By (\ref{eq2.7}) and (\ref{eq2.13}), we have
\begin{equation}\label{eq2.29}
f_k(A)(W_h)\subseteq W_h\quad\mbox{for $A\in \mathcal{A}$ and $k\in S$.}
\end{equation}

By (\ref{eq2.8})  and (\ref{eq2.13}), we have
\begin{equation}\label{eq2.30}
g_k(A)(W_h)\subseteq W_h\quad\mbox{for $A\in \mathcal{A}$ and $k\in S$.}
\end{equation}

Clearly, we have
\begin{equation}\label{eq2.31}
g_h(x)(W_h)\subseteq W_h.
\end{equation}

It follows that $W_h$ is invariant under $f_k(\mathcal{A})$, $g_k(\mathcal{A})$ and $g_h(x)$. Combining a basis of $V^{ann,+}_{2-nd}$ and the linearly independent vectors $u_0$, $u_1$, $\cdots$, $u_{p-1}$, we get a basis of $W_h$. With respect to this basis of $W_h$, the matrix of 
$g_h(\langle x, A \rangle _k)|W_h=[ g_h(x), f_k(A)]|W_h\in End (W_h)$ has the form
$$
 \left(\begin{tabular}{c|c}
$M_{dim\, V^{ann,+}_{2-nd}\times dim\, V^{ann,+}_{2-nd}}$& \mbox{a $(\dim\,V^{ann,+}_{2-nd})\times p$ matrix}\\\hline
$0_{dim\, p\times dim\,V^{ann,+}_{2-nd}}$&$\begin{array}{cccc}
\psi_h (\langle x, A \rangle_k)&\ast &\cdots &\ast\\
0&\psi_h (\langle x, A \rangle_k)&\cdots &\ast \\
\cdots&\cdots&\cdots&\cdots\\
0&0&\cdots&\psi_h (\langle x, A \rangle_k)\end{array}$
\end{tabular} \right)
$$
by (\ref{eq2.8}), where $M_{dim\, V^{ann,+}_{2-nd}\times dim\, V^{ann,+}_{2-nd}}$ is the matrix of $[ g_h(x), f_k(A)]|V^{ann,+}_{2-nd}$ with respect to the basis of 
$V^{ann,+}_{2-nd}$. Hence, we get
$$
tr M=tr\Big([g_h(x), f_k(A)]|V^{ann,+}_{2-nd}\Big)=0
$$
and
$$
0=tr\Big([g_h(x), f_k(A)]|W\Big)=tr M+p\psi_h(\langle x, A \rangle_k)
=p\psi_h(\langle x, A \rangle_k),
$$ 
which implies that
\begin{equation}\label{eq2.32}
\psi_h(\langle x, A \rangle_k)=0 \quad\mbox{for $A\in \mathcal{A}$ and $h$, $k\in S$.}
\end{equation}

Using (\ref{eq2.32}), we have
\begin{eqnarray*}
f_k(A)\Big(g_h(x)(u)\Big)&=&-\Big([g_h(x), f_k(A)] +g_h(x)f_k(A)\Big)(u)\\
&=&-g_h(\langle x, A \rangle_k)(u)+g_h(x)\Big(f_k(A)(u)\Big)\\
&=&-\psi_h(\langle x, A \rangle_k)(u)+g_h(x)\Big(\phi_k (A)(u)\Big)=\phi_k (A)g_h(x)(u),
\end{eqnarray*}
which proves that (\ref{eq2.27}) also holds. By (\ref{eq2.27}), (\ref{eq2.22}) is true.

\hfill\raisebox{1mm}{\framebox[2mm]{}}

\bigskip

The following counterpart of Lie's Theorem claims that if $L$ is a finite dimensional solvable Lie-like algebra$^{2-nd}$ over an algebraically closed field of characteristic $0$, then every finite dimensional non-zero ordinary $L$-module always contains an one dimensional ordinary submodule.

\begin{proposition}\label{pr2.5} {\bf (The Generalization of Lie's Theorem)} Let $L$ be a finite dimensional solvable Lie-like algebra$^{2-nd}$ over an algebraically closed field of characteristic $0$. If $(V, \{f_k\}, \{g_k\}_{k\in S})$  is a finite dimensional ordinary $L$-module, then there exist linear functionals $\phi_k$, $\psi_k: L\to \mathbf{k}$ for $k\in S$ and a nonzero vector $v$ in $V$ such that
\begin{equation}\label{eq36}
f_k(z) (v)=\phi_k (z)v, \quad g_k(z)(v)=\psi_k (z)v \quad\mbox{for $z\in L$ and $k\in S$}. 
\end{equation}
Moreover, either $\psi_k=0$ for all $k\in S$ or $\phi_k= \psi_k $  for all $k\in S$.
\end{proposition}

\medskip
\noindent
{\bf Proof} We use induction on $n:=dim \,L$. If $n=0$, then Proposition~\ref{pr2.5} is clear true.

\medskip
Assume that Proposition~\ref{pr2.5} is true for all finite dimensional solvable Lie-like algebra$^{2-nd}$ of dimensions smaller than $n$. Note that any subspace of $L$ containing
$\mathcal{D}^{2}\mathcal{L}=\sum_{k\in S}\langle L ,\,L\rangle_k$ is an ideal of $L$. Since $n=dim\, L \ge 1$ and $L$ is solvable, 
$L\ne \mathcal{D}^{2}\mathcal{L}$. Therefore, we can find an ideal $\mathcal{A}$ of $L$ such that $\langle L, L\rangle_k\subseteq \mathcal{A}$ for $k\in S$ and $L=\mathcal{A}\oplus \mathbf{k}x$ for some $x\not\in\mathcal{A}$. 

Now $\mathcal{A}$ is a solvable Lie-like algebra$^{2-nd}$ of dimension smaller than $n$ and 
$(V, \{f_k\}, \{g_k\}_{k\in S})$  is a finite dimensional $\mathcal{A}$-module. It follows from induction assumption that
$$
U:=\left\{\, u\in V \, \left |\,\begin{array}{c}\mbox{there exist linear functionals 
$\phi_k' : \mathcal{A}\to \mathbf{k}$}\\\mbox{ and $\psi_k' : \mathcal{A}\to \mathbf{k}$ for $k\in S$ such that}\\\mbox{ $f_k(A)(u)=\phi_k' (A)u$, 
$g_k(A)(u)=\psi_k' (A)u$}\\\mbox{ for $A\in\mathcal{A}$ and $k\in S$}\end{array} \right.\, \right\}\ne 0,
$$

Then we have either $U\bigcap V^{ann,+}_{2-nd}\ne 0$ or $U\bigcap V^{ann,+}_{2-nd}= 0$.

\bigskip
First, we assume that $U\bigcap V^{ann,+}_{2-nd}\ne 0$. By (\ref{eq2.13}) and (\ref{eq2.21}), we have
$$
f_k(x)\Big (U\bigcap V^{ann,+}_{2-nd}\Big)\subseteq U\bigcap V^{ann,+}_{2-nd}
\quad\mbox{for $k\in S$.}
$$
Hence, $\{\, f_k(x)|(U\bigcap V^{ann,+}_{2-nd})\,|\, \mbox{$k\in S$}\,\}$ is a set of commutative linear transformations of the vector space $U\bigcap V^{ann,+}_{2-nd}$ by (\ref{eq1.6}). Thus, there exists 
$0\ne v_0\in U\bigcap V^{ann,+}_{2-nd}$ and $\lambda_k\in \mathbf{k}$ such that
$$
f_k(x)(v_0)=\lambda_k v_0\quad\mbox{for $k\in S$,}
$$
which implies that 
\begin{equation}\label{eq37}
f_k(z)(v_0)=\phi_k(z)v_0\quad\mbox{for $z\in L$ and $k\in S$,}
\end{equation}
where the functional $\phi_k: L\to \mathbf{k}$ is defined by
$$
\phi_k (z)=\left\{\begin{array}{ll}\phi_k' (z)&\mbox{if $z\in \mathcal{A}$,}\\
\lambda _k&\mbox{if $z=x$.}\end{array}\right.
$$

If $g_k(x)(v_0)=0$ for $k\in S$, then
\begin{equation}\label{eq38}
g_k(z)(v_0)=\psi_k(z)v_0\quad\mbox{for $z\in L$ and $k\in S$,}
\end{equation}
where the functional $\psi_k: L\to \mathbf{k}$ is defined by
$$
\psi_k (z)=\left\{\begin{array}{ll}\psi_k' (z)&\mbox{if $z\in \mathcal{A}$,}\\
0&\mbox{if $z\in \mathbf{k}x$.}\end{array}\right.
$$
Hence, (\ref{eq36}) holds for $v:=v_0$ by (\ref{eq37}) and (\ref{eq38}). 

If $g_{h_0}(x)(v_0)\ne 0$ for some $h_0\in S$, then (\ref{eq36}) holds for $v:=g_{h_0}(x)(v_0)$ and $\phi_k=\psi_k=0$ with $k\in S$. In fact, by (\ref{eq1.5}) and (\ref{eq1.6}), we have 
\begin{equation}\label{eq39}
f_k(z)g_{h_0}(x)\big(g_s(y)-f_t(y)\big)=0=g_k(z)g_{h_0}(x)\big(g_s(y)-f_t(y)\big)
\end{equation}
for $z$, $y\in L$ and $k$, $s$, $t\in S$. Since $v_0\in V^{ann,+}_{2-nd}$, $v_0$ is a linear combination of the vectors in the set 
$\{\, \big(g_s(y)-f_t(y)\big)(w)\,|\, \mbox{$s$, $t\in S$, $y\in L$ and $w\in V$}\,\}$. It follows from this fact and (\ref{eq39}) that
$$
f_k(z)\big(g_{h_0}(x)(v_0)\big)=0=g_k(z)\big(g_{h_0}(x)(v_0)\big)
\quad\mbox{for $z\in L$ and $k\in S$,}
$$
which implies that (\ref{eq36}) holds for $v:=g_{h_0}(x)(v_0)$ and $\phi_k=\psi_k=0$ with 
$k\in S$.

\bigskip
We now assume that $U\bigcap V^{ann,+}_{2-nd}= 0$, in which case, we  prove (\ref{eq36}) by cases. Using the same notations as the ones in Proposition~\ref{pr2.3}, we have 
$U=U_{\phi}\cap U_{\psi}$.

\medskip
\underline{\it Case 1:} $f_{h_0}(x)(w)\ne g_{h_0}(x)(w)$ for some $ w\in U$ and some $h_0\in S$, in which case, let
$\tilde{w}:=(f_{h_0}(x)- g_{h_0}(x))(w)$. Then $\tilde{w}\neq 0$. For any $A\in\mathcal{A}$, we have
\begin{eqnarray*}
f_k(A)(\tilde{w})&=&f_k(A)(f_{h_0}(x)-g_{h_0}(x))(w)\\
&=&f_k(A)\Big(f_{h_0}(x)(w)\Big)-f_k(A)\Big(g_{h_0}(x)(w)\Big)\\
&\stackrel{(\ref{eq2.21})}{=}&\phi_k'(A)f_{h_0}(x)(w)-f_k(A)\Big(g_{h_0}(x)(w)\Big)\\
&\stackrel{(\ref{eq2.22})}{=}&\phi_k'(A)f_{h_0}(x)(w)-\phi_k' (A)g_{h_0}(x)(w)\\
&=&\phi_k' (A) (f_{h_0}(x)-g_{h_0}(x))(w)=\phi_k' (A)\tilde{w}
\end{eqnarray*}
or
\begin{equation}\label{eq2.34}
f_k(A)(\tilde{w})=\phi _k'(A)\tilde{w} \quad\mbox{for all $A\in\mathcal{A}$ and $k\in S$.}
\end{equation}

By (\ref{eq1.5}), we have
\begin{equation}\label{eq2.35}
g_k(A)(\tilde{w})=g_k(A)(f_{h_0}(x)- g_{h_0}(x))(w)=0 \quad\mbox{for $A\in\mathcal{A}$ and 
$k \in S$.}
\end{equation}

Let 
$$
\tilde U:=\, \{\, \tilde u\in V\, |\, \mbox{$f_k(A)(\tilde u)=\phi_k' (A)\tilde u$, 
$g_k(A)(\tilde u)=0$ for all $A\in\mathcal{A}$ and $k\in S$} \, \}.
$$
It follows from (\ref{eq2.34}) and (\ref{eq2.35}) that 
$0\ne \tilde{w}\in \tilde U\bigcap V^{ann,+}_{2-nd}$, in which case, we have proved that (\ref{eq36}) holds.

\medskip
\underline{\it Case 2:} $f_h(x)(u)=g_h(x)(u)$ for all 
$h\in S$ and all $u\in U$, in which case, for $A\in\mathcal{A}$, we have
$$ g_k(A)\Big(g_h(x)(u)\Big)\stackrel{(\ref{eq1.5})}{=}
g_k(A)\Big(f_h(x)(u)\Big)\stackrel{(\ref{eq2.21})}{=}\psi_k (A)\Big(f_h(x)(u)\Big)=
\psi_k (A)(g_h(x)(u))$$
or
\begin{equation}\label{eq2.33}
g_h(x)(u)\in U_{\psi} \quad\mbox{for $u\in U=U_{\phi}\cap U_{\psi}$ and $h\in S$.}
\end{equation}

By (\ref{eq2.22}) and (\ref{eq2.33}), $g_h(x)(U)\subseteq U_{\phi}\cap U_{\psi}=U$ for $h\in S$. This fact and (\ref{eq2.21}) prove that $U$ is invariant under both $f_k(x)$ and $g_k(x)$ for $k\in S$. In particular, the set $\{\, f_k(x)|U\,|\, k\in S\,\}$ consisting of commutative linear transformations has a common eigenvector in $U$. Hence,  there exists $\lambda \in\mathbf{k}$ such that
$$ U^{\lambda}: = \{\, u\in U\, |\, \mbox{$f_k(x)(u)=\lambda_k u$ for $k\in S$} \,\}\neq 0.$$
Let $0\ne u_0\in U^{\lambda}$ and $u_m:=g_h(x)^m(u_0)$ for $m\in \mathcal{Z}_{\ge 0}$, where $h$ is a arbitrary fixed element of $S$. We now use induction on $m$ to prove
\begin{equation}\label{eq2.41}
f_k(x)(u_m)\equiv \lambda_k u_m \, mod \, \Big( \displaystyle\sum_{i=0}^{m-1}\mathbf{k}u_i\Big)
\quad\mbox{for $m\in \mathcal{Z}_{\ge 0}$ and $k\in S$.}
\end{equation}

Clearly, (\ref{eq2.41}) holds for $m=0$. Assume that (\ref{eq2.41}) holds for $m$. Then we get
\begin{eqnarray}
f_k(x)(u_{m+1})&=&f_k(x)g_h(x)(u_m)=\Big(-[g_h(x), f_k(x)]+g_h(x)f_k(x)\Big)(u_m)\nonumber\\
\label{eq2.42}&=&-g_h(\langle x, x\rangle_k)(u_m)+g_h(x)f_k(x)(u_m).
\end{eqnarray}

Since $u_m\in U$ for $m\in \mathcal{Z}_{\ge 0}$ and
$\langle x, x\rangle_k \in \mathcal{A}$, we get
\begin{equation}\label{eq2.43}
g_h(\langle x, x\rangle_k)(u_m)=\psi_h '(\langle x, x\rangle_k)u_m .
\end{equation}

Using (\ref{eq2.41}), we get
\begin{eqnarray}
g_h(x)f_k(x)(u_m)&\equiv& g_h(x)(\lambda_k u_m)\, mod \, 
\Big( \displaystyle\sum_{i=0}^{m-1}\mathbf{k}g_h(x)(u_i)\Big)\nonumber\\
\label{eq2.44}&\equiv& \lambda_k g_h(x)(u_m)\, mod \, 
\Big( \displaystyle\sum_{i=1}^m\mathbf{k}u_i\Big)\nonumber\\
&\equiv &\lambda_k u_{m+1}\, mod \, 
\Big( \displaystyle\sum_{i=1}^m\mathbf{k}u_i\Big).
\end{eqnarray}

By (\ref{eq2.42}), (\ref{eq2.43}) and (\ref{eq2.44}), we have
\begin{eqnarray*}
f_k(x)(u_{m+1})&\equiv& -\psi_h '(\langle x, x\rangle_k)u_m +\lambda_k u_{m+1}\, mod \, 
\Big( \displaystyle\sum_{i=1}^m\mathbf{k}u_i\Big)\\
&\equiv& \lambda_k u_{m+1}\, mod \, 
\Big( \displaystyle\sum_{i=1}^m\mathbf{k}u_i\Big),
\end{eqnarray*}
which proves that (\ref{eq2.41}) also holds for $m+1$.

\medskip
Since $V$ is finite dimensional, there exists a positive integer $p$ such that
\begin{equation}\label{eq2.45}
W_h:=\mathbf{k}u_0\oplus \cdots \oplus \mathbf{k}u_{p-1}\quad \mbox{and}\quad 
g_h(x)(W)\subseteq W.
\end{equation}

By (\ref{eq2.41}) and (\ref{eq2.45}), $W_h$ is invariant under both $f_k(x)$ for all $k\in S$ and $g_h(x)$. Since 
$$
[f_k(x), g_h(x)](u_m)\stackrel{(\ref{eq1.4})}{=}
-g_h(\langle x, x\rangle_k)(u_m)\stackrel{(\ref{eq2.43})}{=}-\psi_h '(\langle x, x\rangle_k)u_m ,
$$
the matrix of $[f_k(x), g_h(x)]|W_h$ in the basis $u_0, \cdots , u_{p-1}$ of $W$ is the diagonal matrix:
\[ 
[f_k(x), g_h(x)]|W_h=\left( \begin{array}{cccc}
-\psi_h '(\langle x, x\rangle_k)&0&\cdots &0\\
0&-\psi_h '(\langle x, x\rangle_k)&\cdots &0\\
\cdots&\cdots&\cdots & \\
0&0&\cdots&-\psi_h '(\langle x, x\rangle_k) \end{array}\right)_{p\times p} \]
Hence, $tr\Big([f_k(x), g_h(x)]|W\Big)=-p\psi_h '(\langle x, x\rangle_k)$, which implies that
$\psi_h '(\langle x, x\rangle_k)=0.$ Now for any $u\in U^{\lambda}$, we have
\begin{eqnarray*}
f_k(x)\Big(g_h(x)(u)\Big)&=&\Big([f_k(x), g_h(x)]+g_h(x)f_k(x)\Big)(u)\\
&=&-g_h(\langle x, x\rangle_k)(u)+g_h(x)f_k(x)(u)\\
&=&-\psi_h '(\langle x, x\rangle_k) u+g_h(x)f_k(x)(u)=0+\lambda_k g_h(x)(u),
\end{eqnarray*}
which implies that $g_h(x)(U^{\lambda})\subseteq U^{\lambda}$. This proves that $U^{\lambda}$ is an invariant subspace of $g_h(x)$ for all $h\in S$. Since the set 
$\{\, g_k(x)|U^{\lambda}\,|\, k\in S\,\}$ of commutative linear transformations is commutative, there exists 
$0\ne u'\in U^{\lambda}$ such that $g_k(x)(u')=\mu_k u'$ for $k\in S$ and $\mu_k\in\mathbf{k}$.

Summarizing what we know about $u'\in U^{\lambda}\subseteq U$, we have
$$
f_k(A)(u')=\phi_k ' (A) u', \quad g_k(A)(u')=\psi_k ' (A) u' \quad\mbox{for $A\in\mathbf{A}$ and $k\in S$}
$$
and
$$
f_k(x)(u')=\lambda_k u', \quad g_k(x)(u')=\mu_k u'\quad\mbox{for $A\in\mathbf{A}$ and $k\in S$.}
$$
Hence, if we define linear functionals $\phi_k$, $\psi_k : L\to \mathbf{k}$  by
\begin{eqnarray*}
\phi_k  |\mathcal{A}:=\phi_k ', && \phi_k (x):=\lambda_k,\\
\psi_k |\mathcal{A}:=\psi_k ', && \psi_k (x):=\mu_k
\end{eqnarray*}
for $k\in S$, then
\begin{equation}\label{eq2.46}
f_k(z)(u' )=\phi_k (z) u',\quad g_k(z)(u' )=\psi_k (z)u', \quad\mbox{for $z\in L$ and $k\in S$},
\end{equation}
which implies that (\ref{eq36}) holds for $v:=u'$.

\bigskip
Finally, if $\psi_h\ne 0$ for some $h\in S$, then there exists $y\in L$ such that $\psi_h (y)\ne 0$. Let $v$ be a nonzero vector satisfying (\ref{eq36}). For all $z\in L$ and all $k\in S$, we have
\begin{equation}\label{eq2.47}
g_h(y)f_k(z)(v )=g_h(y)\Big(\phi_k (z) v\Big)=\phi_k (z)\phi_h (y)v
\end{equation}
and
\begin{equation}\label{eq2.48}
g_h(y)g_k(z)(v )=g_h(y)\Big(\psi_k (z) v\Big)=\psi_k (z)\psi_h (y)v
\end{equation}
by (\ref{eq36}). It follows from (\ref{eq2.47}) and (\ref{eq2.48}) that
\begin{equation}\label{eq2.49}
\Big(\phi_k (z)-\psi_k (z)\Big)\psi_h (y)(v )=\Big(g_h(y)f_k(z)-g_h(y)g_k(z)\Big)(v)
\stackrel{(\ref{eq1.5})}{=}0(v)=0.
\end{equation}
Since $\psi_h (y)\ne 0$ and $v\ne 0$, (\ref{eq2.49}) implies that $\phi_k (z)=\psi_k (z)$ for all $z\in L$ and all $k\in S$. It follows from that $\phi_k =\psi_k$ for all $k\in S$.

\hfill\raisebox{1mm}{\framebox[2mm]{}}

\bigskip
As a corollary of Proposition~\ref{pr2.5}, we have the following

\begin{proposition}\label{pr2.6} Let $L$ be a finite dimensional solvable Leibniz algebra over an algebraically closed field of characteristic $0$. If 
$(V, f , g )$ is a finite dimensional nonzero $L$-module, then there exist two
linear functionals $\phi$, $\psi: L\to \mathbf{k}$ and a nonzero vector $v$ in $V$ such that
$$ f(z) (v)=\phi (z)v, \quad g(z)(v)=\psi (z)v \quad\mbox{for $z\in L$}. $$
Moreover, either $(\phi, \psi )=(\phi, \phi )$ or $(\phi, \psi )=(\phi, 0)$.
\end{proposition}

\medskip
\noindent
{\bf Proof} Since a solvable Leibniz algebra is a  solvable Lie-like algebra$^{2-nd}$, 
Proposition~\ref{pr2.6} follows from Proposition~\ref{pr2.5} directly.

\hfill\raisebox{1mm}{\framebox[2mm]{}}

\medskip
Proposition~\ref{pr2.6}, which was announced in \cite{Liu3}, has been used to classify the simple Leibniz algebras with Lie factor $s\ell _2$ in \cite{roots}.

\bigskip
{\bf Acknowledgment} I thank Professor Jean-Louis Loday for his comments on the early version of this paper.

\end{document}